\documentclass[a4paper,12pt]{article}
\usepackage{amssymb}

\newtheorem{theorem}{Theorem}
\newtheorem{lemma}{Lemma}
\newtheorem{remark}{Remark}

\begin{document}

\author{Dumitru B\u{a}leanu\thanks{On leave from The National Institute for Laser, Plasma and Radiation, Physics, Institute of Space Sciences, M\u{a}gurele -- Bucure\c{s}ti, P.O. Box MG-23, R 76911, Romania}\\
\small{\c{C}ankaya University,}\\
\small{Department of Mathematics {\&} Computer Science,}\\
\small{\"{O}gretmenler Cad. 14 06530, Balgat -- Ankara, Turkey}\\
\small{e-mail address: dumitru@cankaya.edu.tr}\\
and\\
Octavian G. Mustafa\\
\small{Department of Mathematics, DAL, University of Craiova, Romania}\\
\small{e-mail address: octawian@yahoo.com}\\
}

\title{On the asymptotic
integration of a class of sublinear fractional differential equations}
\date{}
\maketitle

\noindent{\bf Abstract} We estimate the growth in time of the solutions to a class of nonlinear
fractional differential equations $D_{0+}^{\alpha}\left(x-x_0\right) =f(t,x)$
which includes $D_{0+}^{\alpha}\left(x-x_0\right) =H(t)x^{\lambda}$ with
$\lambda\in(0,1)$ for the case of slowly-decaying coefficients $H$. The proof is based on the triple interpolation
inequality on the real line and the growth estimate reads as $x(t)=o(t^{a\alpha})$
when $t\rightarrow+\infty$ for $1>\alpha>1-a>\lambda>0$. Our result can be thought of as a non--integer
counterpart of the classical Bihari asymptotic integration result for
nonlinear ordinary differential equations. By a carefully designed example we show that in some circumstances such an estimate is optimal.



\section{Introduction}

The fractional calculus deals with the generalization of differentiation
and integration to non-integer orders. Fractional calculus is as old
as the classical differential calculus  and it gained more and more importance during the last
three decades in many fields of science ad engineering
\cite{oldham,miller,Samko,podlubny,Hilfer,zaslavsky,richard,kilbas_book,west,uchaikin,lack}.

The fractional order differential equations play a significant role in modeling the anomalous dynamics of various processes related to complex systems in most areas of science and
engineering. However, only a few steps have been taken toward what
may be called a coherent theory of these equations  in the applied
sciences in a manner analogous to the ordinary case.

We would like to recall that the fractional differential equations have been found effective
to describing various physical phenomena such as damping laws,
rheology, diffusion processes \cite{heymans,mainardi}.
 Various applications of fractional differential equations in the
reaction kinetics of proteins, the anomalous electron transport in
amorphous materials, the dielectrical or mechanical relaxation of
polymers, the modeling of glass--forming liquids and others are presented in
\cite{Metzler,GlNon1995,Sabatier,baleanumagin}. 

On the theoretical side, several fractional
variational principles were formulated and the corresponding
Euler-Lagrange equations consisting of left and the right
derivatives were analyzed during the last few years
\cite{Agrawal1,Atanackovic,Baleanu2008_2,Klimek2002}. The fractional derivatives are also the infinitesimal generators
of a class of translation invariant convolution semigroups which appear universally as attractors \cite{zaslavsky, lack}.
 
New methods for solving the
nonlinear fractional  partial differential equations were developed
recently \cite{momani1,momani2}.

Let us consider the initial value problem for a fractional
differential equation (FIVP) below
\begin{eqnarray}
\left\{
\begin{array}{ll}
D_{0+}^{\alpha}\left(x-x_0\right)(t) =f(t,x(t)),\quad t>0 ,\\
x(0)=x_0,
\end{array}
\right.  \label{fivp}
\end{eqnarray}
where the nonlinearity
$f:D=[0,+\infty)\times\mathbb{R}\rightarrow\mathbb{R}$ is assumed
continuous and $x_0$ is some real number.

The differential operator $D_{0+}^{\alpha}$ in problem (\ref{fivp}) is the
{\it Riemann-Liouville differential operator} of order $0<\alpha<1$,
namely
\begin{eqnarray*}
D_{0+}^{\alpha}(u)(t)=\frac 1{\Gamma (1-\alpha)}\cdot \frac d{dt}\left[
\int\limits_0^t\frac{u(s)}{\left( t-s\right) ^{1-\alpha}}ds\right] ,
\end{eqnarray*}
where $\Gamma(1-\alpha)=\int_{0}^{+\infty}e^{-t}t^{-\alpha}dt$ is the Gamma
function. See \cite[p. 70]{kilbas_book}.

The above fractional differential equation represents a
wide class of dissipative nonlinear nonlocal equations, therefore
the asymptotic behavior of its solutions is an interesting topic of investigation.

The motivation for inserting the initial datum $x_0$ into the differential
operator comes from the physical origin of such mathematical models and
the reader can find comprehensive details in this respect in \cite[p.
80]{podlubny}, \cite[p. 230]{diethelm}. The Riemann-Liouville operator
with an inserted datum is called a {\it Caputo differential operator} --- see
\cite{Caputo}, \cite[p. 91]{kilbas_book}.

Assuming that the FIVP has a solution $x(t)$, the formulas
$\Gamma(\alpha)\Gamma(1-\alpha)=\frac{\pi}{\sin\pi\alpha}$ and
\begin{eqnarray*}
\int_{0}^{t}f(s,x(s))ds=\frac{\sin\pi\alpha}{\pi}
\int_{0}^{t}\frac{1}{(t-s)^{\alpha}}\int_{0}^{s}\frac{f(\tau,x(\tau))}{(s-\tau)^{1-\alpha}}d\tau
ds,
\end{eqnarray*}
see \cite[p. 196]{ahlfors}, allow us to rewrite (\ref{fivp}) via an
integration as
\begin{eqnarray}
\int_{0}^{t}\frac{1}{\Gamma(1-\alpha)(t-s)^{\alpha}}\left[x(s)-x_0-\frac{1}{\Gamma(\alpha)}\int_{0}^{s}\frac{f(\tau,x(\tau))}{(s-\tau)^{1-\alpha}}d\tau\right]ds=0,\label{integral}
\end{eqnarray}
where $t>0$.

We can regard from now on, by means of (\ref{integral}), any solution of
FIVP (\ref{fivp}) as a (continuous) solution of the singular integral
equation
\begin{eqnarray}
x(t)=x_0+\frac{1}{\Gamma(\alpha)}\int_{0}^{t}\frac{f(s,x(s))}{(t-s)^{1-\alpha}}ds,\quad
t>0,\label{integr_ivp}
\end{eqnarray}
and vice versa. This fact is established in a rigorous manner in
\cite{diethelm}. See also the monograph \cite[Theorem 3.24, pp.
199-200]{kilbas_book} (for $\gamma=r=0$ in the original notation).

Given the generality of $f$, we can replace it with $\Gamma(\alpha)^{-1}f$,
which means that we are interested in analyzing the next integral equation
\begin{eqnarray}
x(t)=x_{0}+\int_{0}^{t}\frac{f(s,x(s))}{(t-s)^{1-\alpha}}ds,\quad t\geq0.\label{main}
\end{eqnarray}

Our aim here is to investigate the long--time behavior of the solutions to
a large class of equations (\ref{main}). Noticing that, for a constant $f$
in (\ref{main}), the solution reads as
\begin{eqnarray*}
x(t)=x_0+\frac{f}{\alpha}\cdot t^{\alpha}=O(t^{\alpha})\qquad\mbox{as }t\rightarrow+\infty,
\end{eqnarray*}
we shall compare the growth of solutions to the nonlinear integral
equation (\ref{main}) with the {\it natural} test--function $t^{\alpha}$.

To the best of our knowledge, the literature of fractional differential
equations is scarce with respect to such enterprises. Most of the attempts
regard exact solutions to (very) restricted types of fractional
differential equations. E.g., in \cite{GlNon1995} the linear fractional
differential equation $\tau _0^{-\alpha}D_{0+}^{\alpha}\left[ \Phi (t)\right] +\Phi
(t)-\Phi (0)=0$ has been used to model the complex protein dynamics. The
solutions were expressed using a Mittag--Leffler transcendental function
of order $\alpha$, namely $\Phi (t)=\Phi (0)E_{\alpha}\left( -\left( \frac t{\tau
_0}\right) ^{\alpha}\right)$, see \cite[p. 49]{GlNon1995}. So, $\Phi(t)\sim
t^{-\alpha}$ as the time $t$ increases indefinitely.

The paper is organized as follows. In section \ref{bih_clas} we sketch an asymptotic integration theory for nonlinear ordinary differential equations of integer order that motivates our investigation. In section \ref{prel} we introduce the hypotheses concerning the functional quantity $f(t,x)$ from (\ref{integr_ivp}). The main result is established in section \ref{main_sec}. Section \ref{applications} contains a detailed presentation of the fundamental class of nonlinear differential equations given by $f(t,x)=H(t)x^{\lambda}$, with $\alpha>\lambda>0$. The last section is devoted to an analysis of the sharpness of our asymptotic estimate for the continuable solutions to (\ref{integr_ivp}).

Let us conclude this section by stating that, even though we have chosen to use the simple Riemann-Liouville differential operator inside the nonlinear fractional differential equation (\ref{fivp}) for reasons of clarity in the exposition, the technique presented here works equally well for most of the (far more complicated) differential operators of non-integer order that can be encountered in both the mathematics and physics literature. Another point that has to be emphasized is that, due to the interpolation-type approach described in this paper, we allow the coefficient $H$ to be more flexible (both in the vicinity of the initial moment $0$ and for the large values of $t$) than it is usually allowed in the recent investigations based on the (weighted) H\"{o}lder-type inequalities.

\section{The Bihari asymptotic integration theory}\label{bih_clas}

In two seminal papers \cite{Bihari_1,Bihari_2} on the asymptotic integration of nonlinear ordinary differential equations of integer order, I. Bihari established that, given the equation
\begin{eqnarray*}
x^{\prime\prime}+f(t,x)=0,\quad t\geq t_0>0,
\end{eqnarray*}
where
\begin{eqnarray*}
\vert f(t,x)\vert\leq a(t)\cdot g\left(\frac{\vert x\vert}{t}\right),\qquad x\in\mathbb{R},\thinspace t\geq t_{0},
\end{eqnarray*}
and $f,a,g$ are continuous functions, $g$ is monotone non-decreasing and
\begin{eqnarray*}
\int_{t_0}^{+\infty}a(t)dt<+\infty,
\end{eqnarray*}
all the continuable, real-valued, solutions $x$ satisfy the conditions
\begin{eqnarray*}
\lim\limits_{t\rightarrow+\infty}x^{\prime}(t)=\lim\limits_{t\rightarrow+\infty}\frac{x(t)}{t}=l(x)\in\mathbb{R},
\end{eqnarray*}
that is, they resemble with a certain degree of accuracy to the lines $l(x)\cdot t$.

For the case $l(x)=0$, of interest for applications in mathematical physics, see also the recent contribution \cite{Mustafa2008}.

The core of Bihari's approach relies on using his famous integral inequality, see the presentation from \cite{agarwal_et_al}, to establishing that, if a solution $x$ is globally defined in the future, the quantity
\begin{eqnarray*}
\frac{x(t)}{t},\qquad t\geq t_0,
\end{eqnarray*}
is bounded in $[t_0,+\infty)$.

This implies further that the improper integral
\begin{eqnarray*}
\int_{t_0}^{+\infty}f(s,x(s))ds&=&\lim\limits_{t\rightarrow+\infty}\int_{t_0}^{t}f(s,x(s))ds\\
&=&x^{\prime}(t_0)-\lim\limits_{t\rightarrow+\infty}x^{\prime}(t)
\end{eqnarray*}
is convergent, its value being $x^{\prime}(t_0)-l(x)$.

To search for possible answers to the following question --- are there any
non--integer counterparts of the classical Bihari asymptotic integration
result? ---, we have thus  (1) to produce a quantity like $\frac{x(t)}{t}$ that will ultimately be bounded and further (2) to establish that $\lim\limits_{t\rightarrow+\infty}x^{\prime}(t)$ is finite.

In the following, the quantity from (1) will be
\begin{eqnarray*}
{\cal{B}}(x,t)=\frac{x(t)}{(t+1)^{\alpha}}.
\end{eqnarray*}

The part (2), however, cannot be translated into the non-integer case: this is a consequence of the {\it non-local nature} of the differential operator $D_{0+}^{\alpha}$.

To answer {\it affirmatively} to the preceding question, we have replaced the finiteness of $\lim\limits_{t\rightarrow+\infty}x^{\prime}(t)$, which is equivalent with the finiteness of $\lim\limits_{t\rightarrow+\infty}\frac{x(t)}{t}$, with a condition in the same spirit, namely
\begin{eqnarray*}
\lim\limits_{t\rightarrow+\infty}t^{(1-a)\alpha}\cdot\frac{x(t)}{(t+1)^{\alpha}}\in\mathbb{R}
\end{eqnarray*}
for some $a\in(0,1)$.

\section{Preliminaries}\label{prel}

Several technical requirements for the nonlinearity $f$ will be presented next.

Given $\alpha\in(0,1)$, set $a\in (0,1)$ and $p_2>1$ such that
$\alpha>1-a$ and
\begin{eqnarray}
(1-a)\alpha=\frac{(1-a)(\alpha+{\alpha}^2)}{1+\alpha}<\frac{1}{p_2}<\frac{(1-a)\alpha+{\alpha}^2}{1+\alpha}<\alpha.\label{def_p_2}
\end{eqnarray}

Introduce further $p_1,p_3>1$ via the formulas
\begin{eqnarray*}
\frac{1}{p_1}+\frac{1}{p_2}=\frac{1+(1-a)\alpha}{1+\alpha}\qquad\mbox{and}\qquad
p_3=\frac{1+\alpha}{a\alpha}.
\end{eqnarray*}

Notice that
\begin{eqnarray}
\frac{1}{p_1}+\frac{1}{p_2}+\frac{1}{p_3}=1\label{interpolation},
\end{eqnarray}
which makes $p_1,p_2,p_3$ perfect candidates for the triple interpolation
inequality on the real line, see e.g. \cite[p. 146]{gilbarg}. For further use, remark
that the double inequality regarding $p_2$ from (\ref{def_p_2}) is
equivalent with
\begin{eqnarray}
1-\alpha<\frac{1}{p_1}<\frac{1-(1-a){\alpha}^2}{1+\alpha}.\label{estim_p_1}
\end{eqnarray}

Consider now the continuous, non--decreasing function
$g:[0,+\infty)\rightarrow[0,+\infty)$ with the following features
\begin{eqnarray*}
g(0)=0,\qquad g(\xi)>0\mbox{ when }\xi>0,
\end{eqnarray*}
and --- take $x_0\neq0$ in (\ref{main}) ---
\begin{eqnarray}
\lim\limits_{u\rightarrow+\infty}\frac{W(u)}{u^{a}}=+\infty,\quad
W(u)=\int_{0}^{u}\frac{d\xi}{\left[g\left(\vert
x_0\vert+\xi^{\frac{1}{p_3}}\right)\right]^{p_3}}.\label{g_infinity}
\end{eqnarray}

An example of $g$ is that of $g(u)=u^{1-(1+\zeta)a}$ for some
$\zeta\in(0,1)$ and $a\in\left(0,\frac{1}{1+\zeta}\right)$. It is
easy to see that
\begin{eqnarray*}
W(u)&=&\int_{0}^{u}\frac{d\xi}{\left[\left(\vert
x_0\vert+\xi^{\frac{1}{p_3}}\right)^{1-(1+\zeta)a}\right]^{p_3}}\geq\int_{\vert
x_0\vert^{p_3}}^{u}\frac{d\xi}{\left[\left(2\cdot\xi^{\frac{1}{p_3}}\right)^{1-(1+\zeta)a}\right]^{p_3}}\\
&=&2^{-[1-(1+\zeta)a]p_3}\cdot\frac{u^{(1+\zeta)a}-\vert
x_0\vert^{(1+\zeta)a p_3}}{(1+\zeta)a}\\
&\thicksim& c\cdot
u^{(1+\zeta)a}\quad\mbox{when }u\rightarrow+\infty,
\end{eqnarray*}
where $c=2^{-[1-(1+\zeta)a]p_3}\cdot\frac{1}{(1+\zeta)a}$.

\begin{remark}
\emph{A crucial fact should be emphasized at this point. We shall use in sections \ref{applications}, \ref{fin_sect} the function $g$ given by $g(u)=u^{\lambda}$ for $\alpha>1-a>\lambda$. This restriction upon $\lambda$, which can be recast as $1-\lambda>a$, is natural in the above computation since we can take $(1+\zeta)a=1-\lambda$ for a $\zeta$ small enough.}
\end{remark}

We introduce next the class of Bihari--like nonlinearities $f$ by means of
the inequality
\begin{eqnarray}
\vert f(t,x)\vert\leq h(t)g\left(\frac{\vert
x\vert}{(t+1)^{\alpha}}\right),\label{Bihari_f}
\end{eqnarray}
where the function $h:[0,+\infty)\rightarrow[0,+\infty)$ is continuous and
such that
\begin{eqnarray}
t^{\frac{p_3}{p_1}[1-p_1(1-\alpha)]}\cdot\Vert
h\Vert_{L^{p_2}(0,t)}^{p_3}=O(t^{\alpha})\quad\mbox{when }t\rightarrow+\infty.\label{ground_0}
\end{eqnarray}

\begin{remark}
\emph{The ratio inside $g$ from (\ref{Bihari_f}) is designed, following the tradition of \cite[pp. 277--278]{Bihari_2}, according to our aim of studying the so called {\it long-time behavior} (that is, for $t\rightarrow+\infty$) of solutions to (\ref{fivp}). This is why the choice for the argument of $g$ looks sloppy when $t$ is close to $0$. We emphasize that it is only the quest for clarity and simplicity that has lead us to introducing such a peculiar element $(t+1)^{\alpha}$. A careful inspection of the computations in the following will allow the interested reader to adapt the description of $f$ for the vicinity of $0$. This principle is also valid for the important issue of (any) solution $x$ of (\ref{integr_ivp}) being or not being allowed to take the values $\pm\infty$ for $t=0$. Not to complicate the analysis, we opted for asking $x$ to be continuous as well in $0$. This limitation yields $x(0)=x_{0}$ for all the solutions of (\ref{fivp}), (\ref{integr_ivp}).}
\end{remark}

\begin{remark}
\emph{Our main interest in this paper is about fractional differential equations where the nonlinearity $f(t,x)$ is {\it power-like}. This is motivated by many applied examples: logistic equations used in population dynamics and biochemistry (e.g. polymer growth), Emden-Fowler and Thomas-Fermi type of equations (from stellar dynamics where the chaotic expansions of the outer space gases might be modeled with fractional differential equations). The architecture of the integral $W$ from equation (\ref{g_infinity}) has been chosen, however, as to give a certain freedom to the nonlinearity. One can have even some non-deterministic/non-smooth alterations of the power-like form.}
\end{remark}

The restriction (\ref{ground_0}) can be recast as
\begin{eqnarray}
t^{\frac{p_3}{p_1}[1-p_1(1-\alpha)]}\cdot\left\{\int_{0}^{t}[h(s)]^{p_2}ds\right\}^{\frac{p_3}{p_2}}\leq
M(t+1)^{\alpha},\quad t\geq 0,\label{int_cond_h}
\end{eqnarray}
for some sufficiently large constant $M$.

\begin{remark}
\emph{Notice the inequality of the
exponents of $t$, namely
\begin{eqnarray*}
\frac{p_3}{p_1}[1-p_1(1-\alpha)]<\alpha,
\end{eqnarray*}
which is equivalent with the second part of the double inequality
(\ref{estim_p_1}). This means that} the class of coefficients $h$ is quite
substantial, including, for instance, the space $L^{p_2}(0,+\infty)$. \emph{It is a common practice in the asymptotic analysis of ordinary differential equations of integer order to study carefully the case of} slowly-decaying \emph{coefficients $h$, e.g. $h(t)\thicksim t^{-\varepsilon}$ for some $\varepsilon\in(0,1)$ when $t\rightarrow+\infty$. This issue seems not to be addressed yet in the literature devoted to asymptotic integration of fractional differential equations. For this reason, we hope that our interpolation-type approach will lead to further investigations in this area of research.}
\end{remark}

\section{A Bihari--type result for fractional differential equations}\label{main_sec}

Our main contribution in this work is given in the sequel.
\begin{theorem}\label{theor}
Assume that the function $f$ from (\ref{main}) verifies the conditions
(\ref{Bihari_f}), (\ref{int_cond_h}). Then, all the continuable solutions (that is, continuous solutions defined in $[0,+\infty)$) of the FIVP (\ref{fivp}), (\ref{integr_ivp}) have the asymptotic behavior
\begin{eqnarray}
x(t)=o(t^{a\alpha})\qquad\mbox{as
}t\rightarrow+\infty.
\end{eqnarray}
\end{theorem}

\begin{remark}
\emph{We have refrained here from discussing the delicate matter of the} global existence \emph{of solutions to (\ref{integr_ivp}). This requires the tools of functional analysis and its insertion would have lead us too far afield. However, a long but straightforward adaptation of the procedures from \cite{agarwal_et_al} or \cite{FT} will conclude that such globally defined solutions of the FIVP (\ref{fivp}), (\ref{integr_ivp}) do exist. See also Lemma \ref{lemma_3} of the present paper and \cite{BalMus_conf} for a significant particular case.}
\end{remark}

{\bf Proof.} We split the demonstration into two steps.

In {\it step 1}, we shall establish the raw asymptotic description of solutions to (\ref{fivp}), (\ref{integr_ivp}) given by the formula
\begin{eqnarray}
x(t)=o(t^{\alpha})\quad\mbox{when }t\rightarrow+\infty.\label{first_step_aim}
\end{eqnarray}

To this end, consider $x(t)$ a continuable solution of the integral equation
(\ref{main}) and set $y(t)=x(t)-x_0$ for every $t\geq0$.

We have the estimates
\begin{eqnarray*}
\vert y(t)\vert&\leq&\int_{0}^{t}\frac{1}{(t-s)^{1-\alpha}}\cdot h(s)\cdot
g\left(\frac{\vert y(s)\vert+\vert x_0\vert}{(s+1)^{\alpha}}\right)ds\\
&\leq&\left\{\int_{0}^{t}\left[\frac{1}{(t-s)^{1-\alpha}}\right]^{p_1}\right\}^{\frac{1}{p_1}}\cdot\Vert
h\Vert_{L^{p_2}(0,t)}\cdot[z(t)]^{\frac{1}{p_3}},\quad t>0,
\end{eqnarray*}
where
\begin{eqnarray*}
z(t)=\int_{0}^{t}\left[g\left(\frac{\vert y(s)\vert+\vert
x_0\vert}{(s+1)^{\alpha}}\right)\right]^{p_3}ds.
\end{eqnarray*}

Further --- recall (\ref{int_cond_h}) ---,
\begin{eqnarray}
\vert y(t)\vert^{p_3}&\leq& c\cdot t^{\frac{p_3}{p_1}[1-p_1(1-\alpha)]}\Vert
h\Vert_{L^{p_2}(0,t)}^{p_3}\cdot z(t)\nonumber\\
&\leq&M_1(t+1)^{\alpha}\cdot z(t),\qquad t\geq0,\label{estim_y}
\end{eqnarray}
with $c=[1-p_1(1-\alpha)]^{-\frac{p_3}{p_1}}$ and $M_1=c\cdot M$.

We deduce now that
\begin{eqnarray*}
z^{\prime}(t)\leq [g(\vert y(t)\vert+\vert x_0\vert)]^{p_3}\leq
g\left(\vert x_0\vert+\left[M_1(t+1)^{\alpha}\cdot
z(t)\right]^{\frac{1}{p_3}}\right)^{p_3}.
\end{eqnarray*}

Fix $t_0>0$. Notice that, when $t\in[0,t_0]$, the preceding inequality
reads as
\begin{eqnarray*}
z^{\prime}(t)\leq g\left(\vert x_0\vert+[c_1\cdot
z(t)]^{\frac{1}{p_3}}\right)^{p_3},\quad c_1=M_1(1+t_0)^{\alpha},
\end{eqnarray*}
and
\begin{eqnarray*}
\frac{[c_1z(t)]^{\prime}}{g\left(\vert
x_0\vert+[c_1z(t)]^{\frac{1}{p_3}}\right)^{p_3}}\leq c_1.
\end{eqnarray*}

An integration with respect to $t$ leads to
\begin{eqnarray}
W(c_1z(t))\leq c_1 t\leq c_1 t_0<M_1(1+t_0)^{\alpha+1},\qquad
t\in[0,t_0].\label{intermed_1}
\end{eqnarray}

Replacing $t$ with $t_0$ in (\ref{intermed_1}), we get --- recall
(\ref{estim_y}) ---
\begin{eqnarray}
W(\vert y(t)\vert^{p_3})\leq W(M_1(t+1)^{\alpha}z(t))\leq M_1(t+1)^{\alpha+1},\qquad
t\geq 0.\label{intermed_2}
\end{eqnarray}

The inequality (\ref{intermed_2}) can be rewritten as
\begin{eqnarray}
\frac{W(\vert y(t)\vert^{p_3})}{\left[\vert
y(t)\vert^{p_3}\right]^{a}}\cdot\left[\frac{\vert
y(t)\vert}{(t+1)^{\alpha}}\right]^{a p_3}\cdot\frac{(t+1)^{a\alpha
p_3}}{(t+1)^{\alpha+1}}\leq M_1.\label{intermed_3}
\end{eqnarray}

Notice that the third factor from the second--hand side of
(\ref{intermed_3}) is $1$, that is
\begin{eqnarray}
\frac{W(\vert y(t)\vert^{p_3})}{\left[\vert
y(t)\vert^{p_3}\right]^{a}}\cdot\left[\frac{\vert
y(t)\vert}{(t+1)^{\alpha}}\right]^{a p_3}\leq M_1,\qquad
t\geq0.\label{intermed_4}
\end{eqnarray}

Suppose now, for the sake of contradiction, that $y(t)$ is not $o(t^{\alpha})$
when $t\rightarrow+\infty$. This means that there exist an increasing,
unbounded from above, sequence of positive numbers $(t_n)_{n\geq1}$ and
$\varepsilon>0$ such that
\begin{eqnarray}
\frac{\vert y(t_n)\vert}{(t_{n})^{\alpha}}\geq\varepsilon,\qquad
n\geq1.\label{final_1}
\end{eqnarray}

Since (\ref{final_1}) implies that $\lim\limits_{n\rightarrow+\infty}\vert
y(t_n)\vert=+\infty$, we have --- recall (\ref{g_infinity}) ---
\begin{eqnarray*}
\lim\limits_{n\rightarrow+\infty}\frac{W(\vert
y(t_n)\vert^{p_3})}{\left[\vert
y(t_n)\vert^{p_3}\right]^{a}}=\lim\limits_{u\rightarrow+\infty}\frac{W(u)}{u^{a}}=+\infty.
\end{eqnarray*}

According to the inequality (\ref{intermed_4}), we obtain that
$\lim\limits_{n\rightarrow+\infty}\frac{\vert y(t_n)\vert}{(t_n)^{\alpha}}=0$
and this contradicts (\ref{final_1}).

{\it Step 2.} According to the first step, $y(t)=o(t^{\alpha})$ for all the large values of $t$. This yields, via the L'H\^{o}pital rule, that
\begin{eqnarray*}
\lim\limits_{t\rightarrow+\infty}\frac{z(t)}{t}&=&\lim\limits_{t\rightarrow+\infty}\frac{\int_{0}^{t}\left[g\left(\frac{\vert y(s)\vert+\vert
x_0\vert}{(s+1)^{\alpha}}\right)\right]^{p_3}ds}{t}=\lim\limits_{t\rightarrow+\infty}g\left(\frac{\vert y(t)\vert+\vert
x_0\vert}{(t+1)^{\alpha}}\right)\\
&=&g(0)=0.
\end{eqnarray*}

Further, by means of the estimate (\ref{estim_y}), we conclude that
\begin{eqnarray*}
\vert y(t)\vert=o\left(t^{\frac{1+\alpha}{p_{3}}}\right)=o(t^{a\alpha})\quad\mbox{when }t\rightarrow+\infty.
\end{eqnarray*}

The proof is complete. $\square$

\section{Applications}\label{applications}

To illustrate this work, let us fix $\alpha,\lambda\in(0,1)$ with $\alpha>\lambda$.
There exists $\zeta\in(0,1)$ with the property that
\begin{eqnarray}
\alpha>\frac{\lambda+\zeta}{1+\zeta}.\label{application_1}
\end{eqnarray}
Set also $a=a(\zeta)=\frac{1-\lambda}{1+\zeta}$.

The number $p_2>1$ satisfying (\ref{def_p_2}) can now be introduced.
Noticing that the function
$\zeta\mapsto\frac{\lambda+\zeta}{1+\zeta}=1-a(\zeta)$ is increasing,
the conditions (\ref{def_p_2}) can be simplified
using the stronger restriction
\begin{eqnarray}
\lambda \alpha<\frac{\lambda+\zeta}{1+\zeta}\cdot \alpha<\frac{1}{p_2}&<&\lambda
\alpha+(1-\lambda)\frac{{\alpha}^2}{1+\alpha}=\frac{\lambda
\alpha+{\alpha}^2}{1+\alpha}\label{application_2}\\
&<&\frac{\frac{\lambda+\zeta}{1+\zeta}\cdot \alpha+{\alpha}^2}{1+\alpha}.\nonumber
\end{eqnarray}

In fact, let us fix $p_2$ such that
\begin{eqnarray}
\lambda \alpha<\frac{1}{p_2}<\lambda \alpha+(1-\lambda)\frac{{\alpha}^2}{1+\alpha}.\label{practical_p_2}
\end{eqnarray}
Then, we can replace (\ref{application_1}), (\ref{application_2}) with the
more practical requirements for $\zeta$ given by
\begin{eqnarray*}
\alpha>\frac{\lambda+\zeta}{1+\zeta},\qquad\frac{\lambda+\zeta}{1+\zeta}\cdot
\alpha<\frac{1}{p_2}.
\end{eqnarray*}

Now, the {\it data regarding the numbers $a$, $(p_i)_{1\leq i\leq 3}$} can be summarized in the following algorithm-like structure:

{\it Step 1}\hspace{0.5cm} Take $\alpha$, $\lambda$ such that
\begin{eqnarray*}
1>\alpha>\lambda>0;
\end{eqnarray*}

{\it Step 2}\hspace{0.5cm} Take $p_2$ such that
\begin{eqnarray*}
\left[\lambda\alpha+(1-\lambda)\frac{\alpha^2}{1+\alpha}\right]^{-1}<p_2<(\lambda\alpha)^{-1};
\end{eqnarray*}

{\it Step 3}\hspace{0.5cm} Take $\zeta$ such that
\begin{eqnarray*}
0<\zeta<\min\left\{\frac{\alpha-\lambda}{1-\alpha},\frac{1-\lambda\alpha p_{2}}{\alpha p_{2}-1}\right\};
\end{eqnarray*}

{\it Step 4}\hspace{0.5cm} Set
\begin{eqnarray*}
a=\frac{1-\lambda}{1+\zeta};
\end{eqnarray*}

{\it Step 5}\hspace{0.5cm} Set
\begin{eqnarray*}
p_1=\left[\frac{1+(1-a)\alpha}{1+\alpha}-\frac{1}{p_2}\right]^{-1};
\end{eqnarray*}

{\it Step 6}\hspace{0.5cm} Set
\begin{eqnarray*}
p_3=\frac{1+\alpha}{a\alpha}.
\end{eqnarray*}

Consider the initial value problem for a fractional (general) logistic
equation with sublinear exponents
\begin{eqnarray*}
\left\{
\begin{array}{ll}
D_{0+}^{\alpha}(x-x_0)=H(t)x^{\lambda}\left[1+P(t)x^{-\mu}\right],\qquad t>0,\\
x(0)=x_0,
\end{array}
\right.
\end{eqnarray*}
where $\mu\in(0,\lambda)$, the function
$P:[0,+\infty)\rightarrow\mathbb{R}$ is continuous and bounded, and
\begin{eqnarray*}
H(t)=(t+1)^{-\alpha\lambda}h(t),\qquad h\in L^{p_2}(0,+\infty).
\end{eqnarray*}

Then, all the positive-valued, continuable solutions of the initial value
problem have the asymptotic behavior $x(t)=o(t^{\alpha})$ when
$t\rightarrow+\infty$, according to the first step of Theorem \ref{theor}. This lack of sharpness in the asymptotic estimate is explained in the next paragraph.

Notice that
\begin{eqnarray*}
\left\vert H(t)x^{\lambda}\left[1+P(t)x^{-\mu}\right]\right\vert&\leq&
(t+1)^{-\alpha\lambda}h(t)\cdot(1+\Vert
P\Vert_{\infty})(x^{\lambda}+x^{\eta})\\
&\leq&2(1+\Vert
P\Vert_{\infty})h(t)\left[\frac{1+x}{(t+1)^{\alpha}}\right]^{\lambda}\\
&\leq& h(t)g\left(\frac{x}{(t+1)^{\alpha}}\right),\quad t\geq0,\thinspace x>0,
\end{eqnarray*}
where $\eta=\lambda-\mu\in(0,1)$ and
\begin{eqnarray*}
g(u)=2(1+\Vert P\Vert_{\infty})(1+u)^{\lambda}.
\end{eqnarray*}
Since this function $g$ does not verify the hypothesis $g(0)=0$, we can extract from the inferences of Theorem \ref{theor} only the conclusion of the first step.

Let us observe also that for $a$, $\lambda$ and $p_2$ subjected to the preceding restrictions, that is $0<\lambda<\alpha<1$ and $p_2$ confined by (\ref{practical_p_2}),
all the continuable solutions of the sublinear fractional
differential equation
\begin{eqnarray}
\left\{
\begin{array}{ll}
D_{0+}^{\alpha}(x-x_0)=H(t)x^{\lambda},\qquad t>0,\\
x(0)=x_0,
\end{array}
\right.\label{sublin_cooling}
\end{eqnarray}
where
\begin{eqnarray*}
H(t)=(t+1)^{-\alpha\lambda}h(t),\qquad h\in L^{p_2}(0,+\infty),
\end{eqnarray*}
have the
asymptotic behavior $x(t)=o(t^{a\alpha})$ when $t\rightarrow+\infty$,
according to Theorem \ref{theor}.

As detailed in section
\ref{prel}, the asymptotic integration results for these initial value
problems are manageable illustrations of our technique.
Recently, their superlinear variant, namely the case when $P\equiv0$,
$\alpha\in\left(\frac{1}{2},1\right)$ and $\lambda>1$ has been discussed
in the beautiful contribution \cite{FT}, however, under a {\it quite demanding
restriction} on the coefficient $h$ --- for instance, $h\in
L^{p}(0,+\infty)$ for all $p\in(0,1)$. Such a restriction excludes all hope for an analysis of the long time behavior of solutions in the case of slowly-decaying coefficients. We emphasize the fact that,
even though most real-life applications of (fractional) differential
equations involve nonlinear equations, in the research literature
such enterprises are scarce. In this respect, the contribution
\cite{FT} is quite singular since its methods are mathematically rigorous. On the formal level, a solution of type $O(t^{\alpha})$ as $t\rightarrow+\infty$ for $\alpha=\frac{1}{2}$ and $\lambda=4$ has been presented in \cite[p. 236]{podlubny} for a fractional model of radiative cooling of a semi-infinite body.

\section{A direct analysis}\label{fin_sect}

The final section deals with the crucial issue of {\it sharpness} for our asymptotic estimate in Theorem \ref{theor}. We claim that the asymptotic estimate given in this paper by means of a Bihari-type approach in conjunction with an interpolation inequality is almost optimal. To establish the validity of our claim, we use a method frequently encountered in the asymptotic integration theory of differential equations of integer order, namely we build an "extreme" example and show that, with respect to it, our theorem does the best job ever possible.

Before proceeding with the computations, let us mention that, similarly to the case of the classical Bihari investigation \cite{agarwal_et_al}, we obtain here the following type of error: for an exact solution having the numerical values described by $O(t^{\varepsilon\cdot q})$ when $t\rightarrow+\infty$, where $\varepsilon\in(\varepsilon_{0},1)$ for a prescribed $\varepsilon_{0}>0$ and $q>0$, the asymptotic behavior obtained using inequalities and averaging reads as $o(t^{q})$ for $t\rightarrow+\infty$. Since we can set $\varepsilon_{0}$ as close to $1$ as desired, it is obvious that the estimate provided by such a Bihari-type analysis is optimal.

Consider now the integral equation below
\begin{eqnarray}
x(t)=x_{0}+\int_{0}^{t}\frac{H(s)}{(t-s)^{1-\alpha}}[x(s)]^{\lambda}ds,\quad t>0,\label{analysis_eq}
\end{eqnarray}
where $x_{0}>0$ and $\alpha,\lambda\in(0,1)$ are fixed numbers such that
\begin{eqnarray}
1>\frac{3\alpha}{2}\quad\mbox{and}\quad \alpha>2\lambda.\label{analysis_intermed_0}
\end{eqnarray}
Here, $H:[0,+\infty)\rightarrow[0,+\infty)$ is a continuous function.

The next lemmas address three fundamental issues regarding the investigation of (\ref{analysis_eq}): the sign of $x$, its uniqueness and its global existence in the future.

\begin{lemma}\label{lemma_1}
If $x:[0,T)\rightarrow\mathbb{R}$, where $T\leq+\infty$, is a continuous solution of equation (\ref{analysis_eq}) such that $x(0)=x_{0}$ then
\begin{eqnarray}
x(t)\geq x_{0},\quad t\in[0,T).\label{analysis_intermed_1}
\end{eqnarray}
\end{lemma}

{\bf Proof.} Notice that it is enough to prove that $x$ takes positive values everywhere in $[0,T)$. In fact, if we assume that $x$ has a zero in $[0,T)$ then, given its continuity, there exists the number $t_{0}>0$ such that --- recall that $x(0)=x_{0}>0$ ---
\begin{eqnarray*}
x(t)>0\quad\mbox{for all }t\in[0,t_{0}),\quad x(t_{0})=0.
\end{eqnarray*}

By taking into account (\ref{analysis_eq}), we get
\begin{eqnarray*}
x(t_{0})=x_{0}+\int_{0}^{t_{0}}\frac{H(s)}{(t_{0}-s)^{1-\alpha}}[x(s)]^{\lambda}ds\geq x_{0}.
\end{eqnarray*}

We have reached a contradiction.

Finally, since $x(t)>0$ throughout $[0,T)$, we obtain (\ref{analysis_intermed_1}). $\square$

\begin{lemma}\label{lemma_2}
Given the number $x_{0}>0$, there exists at most one continuous solution $x:[0,T)\rightarrow\mathbb{R}$, where $T\leq+\infty$, of equation (\ref{analysis_eq}) such that $x(0)=x_{0}$.
\end{lemma}

{\bf Proof.} Assume that there exist two continuous solutions $x,y:[0,T)\rightarrow\mathbb{R}$ of equation (\ref{analysis_eq}) with $x(0)=y(0)=x_{0}$. Then, via Lemma \ref{lemma_1}, we have
\begin{eqnarray*}
x(t),y(t)\geq x_{0}\quad\mbox{for all }t\in[0,T).
\end{eqnarray*}

The mean value theorem, applied to the function $k:\left(\frac{x_{0}}{2},+\infty\right)\rightarrow(0,+\infty)$ with the formula $k(\xi)={\xi}^{\lambda}$, yields
\begin{eqnarray*}
\vert x^{\lambda}-y^{\lambda}\vert=\left\vert\frac{\lambda}{{\zeta}^{1-\lambda}}\cdot(\max\{x,y\}-\min\{x,y\})\right\vert
\leq\lambda\left(\frac{2}{x_{0}}\right)^{1-\lambda}\cdot\vert x-y\vert,
\end{eqnarray*}
where $\xi$ lies between $x$ and $y$.

The latter estimate implies that
\begin{eqnarray}
\vert x(t)-y(t)\vert&\leq&\int_{0}^{t}\frac{H(s)}{(t-s)^{1-\alpha}}\cdot\vert [x(s)]^{\lambda}-[y(s)]^{\lambda}\vert ds\nonumber\\
&\leq& c\cdot\int_{0}^{t}\frac{H(s)}{(t-s)^{1-\alpha}}\cdot\vert x(s)-y(s)\vert ds,\quad t\in[0,T),\label{analysis_lipschitz}
\end{eqnarray}
where $c=\lambda(2x_{0}^{-1})^{1-\lambda}$.

The Lipschitz-type formula (\ref{analysis_lipschitz}) allows one to apply the iterative procedure presented e.g. in \cite[pp. 127--131]{podlubny} to conclude that the solutions $x$, $y$ coincide. $\square$

\begin{lemma}\label{lemma_3}
Given the number $x_{0}>0$, there exists at least one continuous solution $x:[0,+\infty)\rightarrow\mathbb{R}$ of equation (\ref{analysis_eq}) such that $x(0)=x_{0}$.
\end{lemma}

{\bf Proof.} The argument is based on a technique developed by the authors in \cite[Theorem 2]{BalMus_conf}.

First, set the numbers $p,q,L>0$ such that
\begin{eqnarray*}
\frac{1}{p}+\frac{1}{q}=1,\quad 1<p<\min\left\{\frac{1}{\alpha},\frac{1}{1-\alpha}\right\},
\end{eqnarray*}
and
\begin{eqnarray*}
L>C(\alpha,p)^q\quad\mbox{for }C(\alpha,p)=p^{\frac{1}{q}-\alpha}\cdot\Gamma(1-p(1-\alpha)).
\end{eqnarray*}

Second, introduce the functions --- here, $f(t,x)=H(t)x^{\lambda}$ ---
\begin{eqnarray*}
k(t)=1+x_0+\int_{0}^{t}\frac{\vert f(s,x_0)\vert}{(t-s)^{1-\alpha}}ds,\quad F(t)=c\cdot H(t),
\end{eqnarray*}
and
\begin{eqnarray*}
H_{L}(t)=k(t)\cdot\exp\left(t+\frac{L}{q}\int_{0}^{t}[k(s)F(s)]^{q}ds\right),\qquad t\geq0,
\end{eqnarray*}
where $c$ has been defined in the preceding proof.

Now, consider the complete metric space ${\cal{X}}=(X,d_L)$, where $X$ is the set of all elements of $C([0,+\infty),[x_0,+\infty))$ that behave as $O(H_{L}(t))$ when $t$ goes to $+\infty$ and
\begin{eqnarray*}
d_L(x,y)=\sup\limits_{t\geq0}\left\{\frac{\vert x(t)-y(t)\vert}{H_{L}(t)}\right\},\qquad x,y\in X.
\end{eqnarray*}

Notice that, as in Lemma \ref{lemma_2}, the following inequality holds:
\begin{eqnarray*}
\vert f(t,x(t))-f(t,y(t))\vert\leq F(t)\vert x(t)-y(t)\vert,\qquad t\geq0,\thinspace x,y\in X.
\end{eqnarray*}

Introduce the operator $T:{\cal{X}}\rightarrow{\cal{X}}$ with the formula
\begin{eqnarray*}
{\cal{T}}(x)(t)=x_0+\int_{0}^{t}\frac{f(s,x(s))}{(t-s)^{1-\alpha}}ds,\qquad x\in X,\thinspace t\geq0.
\end{eqnarray*}

Following verbatim the computations from \cite{BalMus_conf}, we conclude that ${\cal{T}}$ is a contraction of coefficient $C(\alpha,p)\cdot L^{-\frac{1}{q}}$ which implies, via the Banach Contraction Principle, that our problem has a solution globally defined in the future. $\square$

Consider now the particular case of equation (\ref{analysis_eq}) given by
\begin{eqnarray*}
H(t)=(t+1)^{-\alpha\lambda}h(t)=(t+1)^{-\alpha\lambda}\cdot t^{-\frac{1+\varepsilon}{p_{2}}},\quad t\geq1,
\end{eqnarray*}
where $\varepsilon\in(0,1)$ is small enough to have $(1-\lambda)\alpha>\frac{1+\varepsilon}{p_{2}}$. The existence of $\varepsilon$ will be explained in the next paragraph. As regards the definition of $H$, this function is continued downward to $t=0$ in such a way as to keep its sign and continuity in $[0,+\infty)$.

It is supposed that the number $p_{2}$ verifies the double inequality (\ref{practical_p_2}). To see how to choose $\varepsilon$, notice that the inequality below
\begin{eqnarray*}
(1-\lambda)\alpha>\lambda \alpha+(1-\lambda)\frac{{\alpha}^2}{1+\alpha}\quad\left(>\frac{1}{p_{2}}\right)
\end{eqnarray*}
is equivalent with
\begin{eqnarray*}
1>\lambda(2+\alpha).
\end{eqnarray*}
The hypotheses (\ref{analysis_intermed_0}) lead to
\begin{eqnarray*}
1>\frac{3\alpha}{2}>\alpha+\frac{{\alpha}^2}{2}=\frac{\alpha}{2}\cdot(2+\alpha)>\lambda\cdot(2+\alpha),
\end{eqnarray*}
so such a number $\varepsilon$ exists truly.

Let now $x$ be the solution of equation (\ref{analysis_eq}) given by Lemmas \ref{lemma_1}, \ref{lemma_2}, \ref{lemma_3} for $T=+\infty$ --- such a solution exists always ---. By taking into account the formula of $H$, we obtain that
\begin{eqnarray}
x(t)&\geq& x_{0}+\int_{1}^{t}\frac{h(s)}{(t-s)^{1-\alpha}}\cdot\left[\frac{x(s)}{(s+1)^{\alpha}}\right]^{\lambda}ds\nonumber\\
&\geq& x_{0}+\int_{1}^{t}\frac{s^{-\frac{1+\varepsilon}{p_{3}}}}{(t-s)^{1-\alpha}}\cdot\left[\frac{x_{0}}{(s+1)^{\alpha}}\right]^{\lambda}ds\nonumber\\
&\geq&\int_{1}^{t}\frac{ds}{(t-s)^{1-\alpha}}\cdot\frac{1}{t^{\frac{1+\varepsilon}{p_{3}}}}\cdot\frac{x_{0}^{\lambda}}{(t+1)^{\alpha\lambda}}\nonumber\\
&\thicksim& c\cdot t^{(1-\lambda)\alpha-\frac{1+\varepsilon}{p_{3}}}\quad\mbox{when }t\rightarrow+\infty,\label{analysis_intermed_expo}
\end{eqnarray}
where $c=\frac{x_{0}^{\lambda}}{\alpha}$.

According to our theorem, $x(t)=o(t^{\alpha a})$ for $t\rightarrow+\infty$, where $\lambda<1-a<\alpha$. Since we can fix $a$ in this range freely, assume that $a$ is taken such that
\begin{eqnarray}
\alpha>2\lambda>1-a>\lambda.\label{analysis_intermed_2}
\end{eqnarray}
This choice is in perfect agreement with the hypotheses (\ref{analysis_intermed_0}).

The restriction (\ref{analysis_intermed_2}) leads to --- recall (\ref{practical_p_2}) ---
\begin{eqnarray*}
\frac{1}{p_{2}}>\lambda \alpha>(1-\lambda-a)\alpha.
\end{eqnarray*}
Further, we have
\begin{eqnarray*}
a\alpha>(1-\lambda)\alpha-\frac{1}{p_{2}}>(1-\lambda)\alpha-\frac{1+\varepsilon}{p_{2}}>0.
\end{eqnarray*}

We would like to evaluate the difference between the exponent $a\alpha$ of $t$ in the asymptotic formula of the solution $x$ given by our theorem and the exponent $(1-\lambda)\alpha-\frac{1+\varepsilon}{p_{2}}$ provided by (\ref{analysis_intermed_expo}). This reads as
\begin{eqnarray*}
&&a\alpha-\left[(1-\lambda)\alpha-\frac{1+\varepsilon}{p_{2}}\right]=\frac{1+\varepsilon}{p_{2}}+(a+\lambda-1)\alpha\\
&&<(1+\varepsilon)\left[\lambda \alpha+(1-\lambda)\frac{{\alpha}^{2}}{1+\alpha}\right]+(a+\lambda-1)\alpha\\
&&=\alpha\left\{[a+(2+\varepsilon)\lambda-1]+\frac{1+\varepsilon}{1+\alpha}[\alpha(1-\lambda)]\right\}\\
&&<\alpha\{[a+(2+\varepsilon)\lambda-1]+(1+\varepsilon)\alpha\}=\Omega \alpha.\\
\end{eqnarray*}

Fix now $\eta_{0}\in(0,1)$ and replace the first of hypotheses (\ref{analysis_intermed_0}) with a stronger condition, namely
\begin{eqnarray*}
\eta_{0}>\frac{7\alpha}{2}.
\end{eqnarray*}
Thus, we have
\begin{eqnarray*}
\Omega<\eta_{0}a.
\end{eqnarray*}
In fact, this inequality is equivalent with
\begin{eqnarray*}
(1-\eta_{0})a+(2+\varepsilon)\lambda+(1+\varepsilon)\alpha<1.
\end{eqnarray*}

The latter inequality is valid as a consequence of the estimate below --- recall the second of hypotheses (\ref{analysis_intermed_0}) ---
\begin{eqnarray*}
&&(1-\eta_{0})a+(2+\varepsilon)\lambda+(1+\varepsilon)\alpha<(1-\eta_{0})+3\lambda+2\alpha\\
&&<1-\eta_{0}+\frac{3\alpha}{2}+2\alpha<1.
\end{eqnarray*}

The conclusion of this analysis of equation (\ref{analysis_eq}) is that the solution $x$, if computed numerically, will be at least as big as $O(t^{(1-\eta)\cdot a \alpha})$ for $t\rightarrow+\infty$, where $\eta\in(0,\eta_{0})$ and $\eta_{0}$ is a {\it prescribed} quantity, while the description of $x$ given by our theorem reads as $o(t^{a\alpha})$ when $t\rightarrow+\infty$. In other words, the theorem provides a sharp asymptotic estimate for the solutions of (\ref{fivp}), (\ref{integr_ivp}) and, in some particular cases, this estimate is optimal.

\section{Acknowledgments}
The work on this paper has been done during the visit of OGM to
\c{C}ankaya University in the fall of 2008. OGM is grateful to the people
from the Department of Mathematics and Computer Science for the friendly
and enthusiastic working atmosphere. Both authors are grateful to the referees for their helpful comments.

\end{document}